\documentclass[12pt]{amsart}
\usepackage[applemac]{inputenc}
\usepackage{bm}

\usepackage[dvipsnames,svgnames,x11names,hyperref]{xcolor}
\usepackage[hyphens]{url}
\usepackage[colorlinks=true,linkcolor=Maroon,citecolor=blue,urlcolor=blue,hypertexnames=false,linktocpage]{hyperref}
\usepackage{bookmark}
\usepackage{amsmath,thmtools,mathtools}
\mathtoolsset{showonlyrefs=true}

\usepackage{blindtext}

\newcounter{mgncount}

\usepackage{fancyhdr}
\usepackage{esint}
\usepackage{enumerate}

\usepackage{pictexwd,dcpic}
\usepackage{graphicx}
\usepackage{caption}
\usepackage{graphicx}

\newtheorem*{thm-}{Theorem}
\declaretheorem[name=Theorem,numberwithin=section]{thm}

\declaretheorem[name=Lemma,sibling=thm]{lemma}

\numberwithin{equation}{section}

\newcommand{\Rn}{\mathbb{R}^{n+1}}
\newcommand{\Sn}{\mathbb{S}^n}

\newcommand{\pf}[1]{\begin{proof}#1 \end{proof}}
\newcommand{\eq}[1]{\begin{equation}\begin{alignedat}{2} #1 \end{alignedat}\end{equation}}

\newcommand{\ra}{\rightarrow}

\makeatletter
\@namedef{subjclassname@2020}{%
	\textup{2020} Mathematics Subject Classification}
\makeatother

\begin{document}
	\title[isotropic $L_{p}$ dual Minkowski problem]
	{On the uniqueness of solutions to the isotropic $L_{p}$ dual Minkowski problem}
	\author[Y. Hu, M. N. Ivaki]{Yingxiang Hu, Mohammad N. Ivaki}
	\maketitle
\begin{abstract}
We prove that the unit sphere is the only smooth, strictly convex solution to the isotropic $L_p$ dual Minkowski problem
	\begin{align*}
		h^{p-1} |D h|^{n+1-q}\mathcal{K}=1,
	\end{align*} 
provided $(p,q)\in (-n-1,-1]\times [n,n+1)$.
\end{abstract}
	
\section{Introduction}
An important question in convex geometry is the uniqueness or the non-uniqueness of the origin-centred spheres as solutions to the isotropic $L_{p}$ dual Minkowski problem
	\begin{align}\label{Lpq-Minkowski-problem}
		h^{p-1} |D h|^{n+1-q} \mathcal{K}=c, \quad c\in (0,\infty).
	\end{align}
The $L_p$ dual Minkowski problem was first introduced by Lutwak, Yang and Zhang \cite{LYZ18}, acting as a bridge which connects the $L_p$-Minkowski problem to the dual Minkowski problem. The former, the $L_p$-Minkowski problem, was introduced by Lutwak in his influential paper \cite{Lut93} three decades ago, and since then has been extensively investigated; e.g., \cite{LO95,LYZ04,Sta02,Sta03,CW06,BLYZ12,BLYZ13,JLW15,Zh15,HLW16,BIS19,BBCY19,Li19,CHLL20,Mil21,KM22,IM23b} and \cite{Sch14}. The latter, the dual Minkowski problem, was proposed recently by Huang et. al in \cite{HLYZ16}, and further studied in \cite{Zh17,Zh18,BHP18,CL18,HP18,HJ19,LSW20}. There has been significant progress on the $L_p$ dual Minkowski problem after the paper \cite{LYZ18} such as \cite{HZ18,BF19,CHZ19,CL21,LLL22}; however, the complete answer to the uniqueness and non-uniqueness question, as stated above, has been elusive in the most interesting case: without the \emph{origin-symmetry} assumption.
	
Here are the known uniqueness and non-uniqueness results for the isotropic $L_{p}$ dual Minkowski problem:
\begin{itemize}
		\item \cite{BCD17}, uniqueness of solutions for $-(n+1)\leq p< 1$ and $q=n+1$ (see also \cite{And99,And03,Sar22});
		\item \cite{HZ18}, uniqueness of solutions for $p>q$;
		\item \cite{CHZ19}, uniqueness of origin-symmetric solutions for 
		\[
			-(n+1)\leq p < q\leq \min\{n+1,n+1+p\};
		\]
		\item \cite{CL21}, uniqueness of solutions for $1<p<q\leq n+1$, or $-(n+1)\leq p<q<-1$, or the uniqueness of solutions up to rescaling for $p=q$;
		\item \cite{LW22}, complete classification for $n=1$;
		\item \cite{CCL21}, non-uniqueness of solutions under any of the following assumptions:
		\begin{enumerate}[(i)]
			\item $q-2(n+1) >p\geq 0$;
			\item $q>0$ and $-q^\ast <p <\min\{ 0,q-2n-2\}$, where 
			\[
				q^\ast:=\left\{\begin{aligned}&\frac{q}{q-n}, \quad \text{if $q\geq n+1$}\\
				&\frac{nq}{q-1}, \quad \text{if $1<q<n+1$}\\
				&+\infty, \quad \text{if $0<q\leq 1;$}
			\end{aligned}\right.
			\]
			\item $p+2(n+1) <q\leq 0$.
		\end{enumerate}
\end{itemize}
	
In the recent work \cite{IM23a}, employing the local Brunn-Minkowski inequality, the following uniqueness was proved.
\begin{thm-}
Let $n\geq 2$ and assume $-(n+1)\leq p$ and $q\leq n+1$, with at least one being strict. Suppose $\mathcal{M}^n$ is a smooth, strictly convex, origin-centred hypersurface such that $h^{p-1}|Dh|^{n+1-q}\mathcal{K}=c$ with $c>0$. Then $\mathcal{M}^n$ is an origin-centred sphere.
\end{thm-}
	
Here, we also employ the local Brunn-Minkowski inequality as our main tool to establish the following uniqueness result.

\begin{thm}\label{s1:cor-Lp-dual-Minkowski}
Let $n\geq 2$. Suppose $\mathcal{M}^n$ is a smooth, strictly convex hypersurface with $h>0$, such that $h^{p-1}|Dh|^{n+1-q}\mathcal{K}=1$. Suppose either
	\begin{enumerate}
		\item $-(n+1)< p\leq -1$ and $n\leq q\leq n+1$,
		\item or $-(n+1)\leq p\leq -n$ and $1\leq q< n+1$.
	\end{enumerate}
Then $\mathcal{M}^n$ is the unit sphere.
\end{thm}
	
\section{Background}
\subsection{Convex geometry}
Let $(\mathbb R^{n+1},\delta:=\langle\, ,\rangle,D)$ denote the Euclidean space with its standard inner product and flat connection, and let $(\mathbb S^n,\bar{g},\bar{\nabla})$ denote the unit sphere equipped with its standard round metric and Levi-Civita connection. 
	
Suppose $K$ is a smooth, strictly convex body in $\mathbb R^{n+1}$ with the origin in its interior. Write $\mathcal{M}=\mathcal{M}^n=\partial K$ for the boundary of $K$. The Gauss map of $\mathcal{M}$, denoted by $\nu$, takes the point $p\in \mathcal{M}$ to its unique unit outward normal $x=\nu(p)\in \mathbb S^n$. The support function of $K$ is defined by 
	\begin{align*}
		h(x):=\max\{\langle x,y\rangle:~y\in K\}, \quad x\in \mathbb S^n.
	\end{align*}
The inverse Gauss map $X=\nu^{-1}:\mathbb S^n\ra \mathcal{M}$ is given by
	\begin{align*}
		X(x)=Dh(x)=\bar{\nabla} h(x)+h(x)x, \quad x\in \mathbb S^{n}.
	\end{align*}
The support function can also be expressed as
	\begin{align*}
		h(x)=\langle X(x),x\rangle=\langle \nu^{-1}(x),x\rangle, \quad x\in \mathbb S^n.
	\end{align*}
The radial function of $K$ is defined by
	\begin{align*}
		r(x):=|X(x)|=(|\bar{\nabla} h(x)|^2+h^2(x))^\frac{1}{2}.
	\end{align*}
Moreover, the Gauss curvature of $\mathcal{M}$ is defined by
	\begin{align*}
		\frac{1}{\mathcal{K}(x)}:=\left.\frac{\det(\bar{\nabla}^2 h+\bar{g}h)}{\det(\bar{g})}\right|_x, \quad x\in \mathbb S^n.
	\end{align*}
	
Note that the matrix $A[h]:=\bar{\nabla}^2 h+h\bar{g}=D^2 h|_{T\mathbb S^n}$ is positive-definite. The eigenvalues of the matrix $A[h]$ with respect to the metric $\bar{g}$, denoted by $\lambda=(\lambda_1,\ldots,\lambda_n)$, are the principal radii of curvature at the point $X(x)\in \mathcal{M}$. 
Then $\sigma_n=\mathcal{K}^{-1}=\Pi_i \lambda_i$. The curvature equation \eqref{Lpq-Minkowski-problem} can be reformulated as the following Monge-Amp\'ere equation:
	\eq{
		h^{1-p} |D h|^{q-n-1} \det(\bar{\nabla}^2 h+\bar{g} h)=c.
	}
The polar body of $K$ is defined by 
	\begin{align*}
		K^\ast:=\{y\in \mathbb R^{n+1}: \langle x,y\rangle\leq 1\, ~\forall x\in K\}.
	\end{align*}
It is well-known that $K^\ast$ is also a smooth, strictly convex body $K^\ast$ in $\mathbb R^{n+1}$ with the origin in its interior. Moreover, the following identity holds
	\begin{equation}\label{polar-dual}
		\frac{h^{n+2}(x) (h^\ast(x^{\ast}))^{n+2}}{\mathcal{K}(x)\mathcal{K}^\ast(x^\ast)}=1.
	\end{equation}
Here $h^{\ast}$ and $\mathcal{K}^{\ast}$ denote respectively the support function and Gauss curvature of $K^{\ast}$, and $x^{\ast}:=X(x)/|X(x)|$.
	
Finally, let us introduce the measure $dV:=h\sigma_nd\mu$,
where $\mu$ is the spherical Lebesgue measure of the unit sphere $\mathbb S^n$. Then the measure $\sigma_n d\mu$ is the surface-area measure of $K$, and $dV=h\sigma_n d\mu$ is a constant multiple of the cone-volume measure of $K$. We refer to \cite{Sch14} for additional background.

\subsection{Centro-affine geometry}
In this section, we recall some basics from centro-affine geometry.  For the related concepts, we refer the reader to \cite{LHSZ15,NS94} and, in particular, to the excellent paper by Milman \cite{Mil21}.

Let $X: \Sn\to \mathcal{M}$ be a smooth embedding of $\mathcal{M}$ (which we consider it to be $Dh$ as in the previous section), and consider the transversal normal field $\xi(x):=X(x)$ (the centro-affine normal). The transversal vector $\xi$ induces the volume form $V$ (as in the previous section), a connection $\nabla$, as well as a metric $g^{\xi}$ on $\Sn$ as follows:
\begin{align*}
V(e_1,\ldots, e_{n})=\det (dX(e_1),\ldots, dX(e_{n}),\xi),\quad e_i\in T\Sn,	
\end{align*}
\eq{\label{nabla structure}
D_udX(v)=dX(\nabla_uv)-g^{\xi}(u,v)\xi,\quad u,v\in T\Sn.
}
Note that $g^{\xi}$ is symmetric and positive-definite. Moreover, while $\nabla$ is not the Levi-Civita connection of $g$, it is torsion-free and
\eq{\label{integration by parts}
\nabla V\equiv 0.
}

The conormal field $\xi^\ast: \Sn\to (\Rn)^{\ast}\sim \Rn$ is the unique smooth vector field to the dual space of $\Rn$, such that $\langle \xi^\ast,dX\rangle=0$ and $\langle \xi, \xi^\ast\rangle=1$. Moreover, $\xi^\ast$ is an immersion and transversal to its image, and it induces a bilinear form and a torsion-free connection on $\Sn,$ 
\eq{
D_ud\xi^\ast(v)=d\xi^{\ast}(\nabla_u^{\ast}v)-g^{\xi^\ast}(u,v)\xi^\ast ,\quad u,v\in T\Sn.
}
We furnish all geometric quantities associated with $\xi^{\ast}$ with $\ast$.

It is known that $g^{\xi}=g^{\xi^{\ast}}$ and that the two connections $\nabla^{\ast}$ and $\nabla$ are conjugate with respect to $g^{\xi}$:
\eq{
ug^{\xi}(v_1,v_2)=g^{\xi}(\nabla_uv_1,v_2)+g^{\xi}(v_1,\nabla^{\ast}_uv_2) \quad u,v_1,v_2\in T\Sn.
}
Moreover,
by \cite[Proposition 4.2]{Mil21} (or taking the inner production of \eqref{nabla structure} with $\nu$), we find
\eq{
g^{\xi}=g^{\xi^{\ast}}=\frac{A[h]}{h}:=g.
}

For a smooth function $f:\Sn\to \mathbb{R}$, the Hessian and Laplacian with respect to $(\nabla,g)$ are defined as
\[\operatorname{Hess}f(u,v)=\nabla df(u,v)=v(uf)-df(\nabla_vu)\]
and
$
\Delta f=\operatorname{div}_g(\nabla f)=\sum_{i} g(\nabla_{e_i}\nabla f,e_i),
$
where $\{e_i\}_{i=1}^n$ is a local $g$-orthonormal frame of $T\Sn$. 

We write $\operatorname{Hess}^\ast$ and $\Delta^{\ast}$ respectively for the Hessian and Laplacian with respect to $(\nabla^{\ast},g)$.
Since $\nabla,\nabla^{\ast}$ are conjugate, we have
\begin{align*}
v(uf)=vg(\nabla f,u)=g(\nabla_v\nabla f,u)	+df(\nabla^{\ast}_vu).
\end{align*}
Therefore, we obtain
\eq{
\Delta f=\operatorname{tr}_g\operatorname{Hess}^{\ast}f,\quad \Delta^{\ast} f=\operatorname{tr}_g\operatorname{Hess}f.
}
By \cite[Proposition 4.2]{Mil21}, we have
\eq{
\operatorname{Hess}^{\ast}f+gf=\frac{1}{h}\left(\bar{\nabla}^2(hf)+\bar{g}hf\right)=\frac{A[hf]}{h}.
}
Let us define
\eq{
Q(u,v)=\nabla^{\ast}_vu-\nabla_v u \quad \forall u,v \in T\Sn.
}
Then by \cite[(6.2)]{LHSZ15},
\eq{
\operatorname{tr}_{g}Q=-\nabla\log\left(\frac{h^{n+2}}{\mathcal{K}}\right).
}
In particular, we have
\eq{\label{change of laplacian}
(\Delta- \Delta^{\ast})f=-\sum_iQ(e_i,e_i)f=d\log \frac{h^{n+2}}{\mathcal{K}}(\nabla f).
}

We conclude this section by recalling the local Brunn-Minkowski inequality, reformulated in the language of centro-affine geometry (cf.  \cite{Mil21}):
Let $f\in C^1(\Sn)$. Then
\eq{\label{local BM}
n\int f^2 dV\leq \int |\nabla f|_g^2 dV+n\frac{(\int fdV)^2}{\int dV}.
}
The equality holds if and only if for some $w\in \Rn$
\eq{
f(x)=\langle \frac{x}{h(x)},w\rangle\quad \forall x\in \Sn.
}
Moreover, by \cite[(5.9)]{Mil21} we also have 
\eq{\label{local BM 2}
n\int |\nabla f|_g^2 dV\leq \int (\Delta f)^2 dV \quad \forall f\in C^2(\Sn).
}
\section{Uniqueness}	
The following identity is at the heart of our approach to employing the local Brunn-Minkowski inequality.
\begin{thm}\label{main identity} There holds
	\eq{\label{main inq}
	\Delta X+nX=h\bar{\nabla}\log \frac{h^{n+2}}{\mathcal{K}}.
	}
In particular,
\eq{
n\int X dV= \int h\bar{\nabla}\log \frac{h^{n+2}}{\mathcal{K}} dV.
}
\end{thm}
\pf{ Let $w\in \mathbb{R}^{n+1}$ be a fixed vector.
By the centro-affine Gauss equation for $\xi=X$ (cf. \cite[Section 3.8]{Mil21}), we have
\eq{
\Delta^{\ast}\langle X,w\rangle+n\langle X,w\rangle=0.
}
Now let $\{v_i\}_{i=1}^n$ be a local orthonormal frame of $T\mathbb{S}^n$ that diagonalizes $A[h]$ at $x_0$ and $A[h]|_{x_0}(v_i,v_j)=\delta_{ij}\lambda_i$. Define $e_i=\sqrt{\frac{h}{\lambda_i}}v_i$, $i=1,\ldots, n$. Then we have $g|_{x_0}(e_i,e_j)=\delta_{ij}$.
Hence, by \eqref{change of laplacian}, at $x_0$ we have
\eq{
\Delta \langle X,w\rangle+n\langle X,w\rangle&=(\Delta -\Delta^{\ast}) \langle X,w\rangle\\
&= g( \nabla\log \frac{h^{n+2}}{\mathcal{K}},\nabla  \langle X,w\rangle)\\
&= g( \nabla\log \frac{h^{n+2}}{\mathcal{K}},\lambda_i  \langle e_i,w\rangle e_i)\\
&=\lambda_i  \langle e_i,w\rangle  d\log \frac{h^{n+2}}{\mathcal{K}}(e_i)\\
&= \langle  \bar{\nabla}\log \frac{h^{n+2}}{\mathcal{K}},hw\rangle. 
}
The second identity follows from integrating \eqref{main inq} against $dV$.
}
\begin{lemma}\label{ineq local BM - a}
Let $0<f\in C^2(\mathbb{S}^n)$. Then
\eq{
\int  f^2\left(\langle \bar{\nabla}\log \frac{h^{n+2}}{\mathcal{K}}, hX\rangle-|X|^2|\nabla\log f|_g^2\right) dV\leq n\frac{|\int  fXdV|^2}{\int dV}.
}	
\end{lemma}
\pf{Let $\{E_k\}_{k=1}^{n+1}$ be an orthonormal basis of $\mathbb{R}^{n+1}$. We define 
\eq{
f_k=f\langle X,E_k\rangle\quad k=1,\ldots,n+1.
}
In view of \autoref{main identity}, we have
\begin{align*}
\Delta f_k+nf_k=&f\langle \bar{\nabla}\log \frac{h^{n+2}}{\mathcal{K}}, hE_k\rangle+\langle X,E_k\rangle\Delta f+2g(\nabla f,\nabla \langle X,E_k\rangle).
\end{align*}
Therefore,
\begin{align}\label{fk local BM}
\sum_kf_k(\Delta f_k+nf_k)=&f^2\langle \bar{\nabla}\log \frac{h^{n+2}}{\mathcal{K}}, hX\rangle+f|X|^2\Delta f\\
&+fg(\nabla f,\nabla |X|^2).
\end{align}
Moreover, by integration by parts (cf. \eqref{integration by parts}), there holds
\eq{\label{integration by parts-}
\int |X|^2f\Delta f+fg(\nabla f, \nabla |X|^2)dV=-\int |X|^2|\nabla f|_g^2 dV.
}
By the local Brunn-Minkowski inequality (see \eqref{local BM}), we have
\eq{
\sum_k\int f_k(\Delta f_k+nf_k)dV\leq n \sum_k\frac{\langle\int fXdV,E_k\rangle^2 }{\int dV}.
}
Thus the claim follows from \eqref{fk local BM} and \eqref{integration by parts-}.
}
\begin{lemma} \label{ineq local BM - b}
	Suppose $\varphi: (0,\infty)\to (0,\infty)$ is $C^1$-smooth and $f=\varphi(r)$. Then we have
	\eq{
\int  f^2\langle \bar{\nabla}\log \frac{h^{n+2}}{\mathcal{K}}-(r(\log \varphi)')^2\bar{\nabla}\log r, hX\rangle dV\leq n\frac{|\int  fXdV|^2}{\int dV}.
}
\end{lemma}
\pf{Let $\{v_i\}_{i=1}^n$ and $\{e_i\}_{i=1}^n$ be as in the proof of \autoref{main identity}. We calculate
\eq{
e_i(\log f)=(\log \varphi)'e_ir= \frac{(\log \varphi)'}{r}\lambda_i\langle e_i,X\rangle=\frac{(\log \varphi)'}{r}\sqrt{h\lambda_i}\langle v_i,X\rangle,
}
and
\eq{
r^2|\nabla \log f|_g^2=((\log \varphi)')^2h\lambda_i(v_ih)^2 = (r(\log \varphi)')^2\langle \bar{\nabla}\log r, hX\rangle.
}
Now the inequality follows from \autoref{ineq local BM - a}.
}

\pf{[Proof of \autoref{s1:cor-Lp-dual-Minkowski}] 
Let $\alpha=q-n-1$. Due to \autoref{ineq local BM - b} with $\varphi(r)=r^{q-n-1}$, and our assumption
$
h^{n+2}\mathcal{K}^{-1}=h^{n+1+p}r^{n+1-q},
$ 
we obtain
\eq{
&(n+1+p)\int  r^{2\alpha}|\bar{\nabla}h|^2dV\\
\leq~& \alpha(\alpha+1)\int  r^{2\alpha} \langle \bar{\nabla}\log r,h\bar{\nabla}h\rangle dV+ n\frac{|\int  r^{\alpha}XdV|^2}{\int  dV}.
}
Assuming $\alpha^2+\alpha\leq0$ (i.e. $n\leq q\leq n+1$) we obtain
\eq{
(n+1+p)\int  r^{2\alpha}|\bar{\nabla}h|^2dV\leq n\frac{|\int  r^{\alpha}XdV|^2}{\int  dV}.
}
Moreover, by using $\bar{\Delta} x+nx=0,$
\eq{
\int  r^{\alpha}XdV=\int Xh^{p}d\mu=\frac{n+1+p}{n}\int  r^{\alpha}\bar{\nabla}h dV.
}
Hence, due to $n+1+p>0$,
\eq{
\int  r^{2\alpha}|\bar{\nabla}h|^2dV\leq \frac{n+1+p}{n}\frac{|\int  r^{\alpha}\bar{\nabla}hdV|}{\int  dV}.
}
We may rewrite this inequality as
\eq{
\int  \left|r^{\alpha}\bar{\nabla}h-\frac{\int r^{\alpha}\bar{\nabla}h dV}{\int  dV}\right|^2dV\leq \frac{p+1}{n}\frac{|\int r^{\alpha}\bar{\nabla}h dV|^2}{\int  dV}.
}
Thus $h$ is constant, provided $-(n+1)<p\leq -1$ and $n\leq q \leq n+1$. 

In view of \eqref{polar-dual}, the polar body $K^\ast$ satisfies the following isotropic $L_{-q}$ dual Minkowski problem:
	\begin{align*}
		(h^{\ast})^{-1-q}|Dh^\ast|^{n+1+p}\mathcal{K}^\ast=1.
	\end{align*}
Hence, the uniqueness result also holds when $n\leq -p\leq (n+1)$ and $-(n+1)< -q\leq -1$.
	}

\section*{Acknowledgment}
 The first author's work was supported by the National Key Research and Development Program of China 2021YFA1001800 and the National Natural Science Foundation of China 12101027. Both authors were supported by the Austrian Science Fund (FWF): Project P36545. 
	
\providecommand{\bysame}{\leavevmode\hbox to3em{\hrulefill}\thinspace}
\providecommand{\href}[2]{#2}

	\vspace{10mm}
	\textsc{School of Mathematical Sciences, Beihang University, Beijing 100191, China}
	\email{\href{mailto:huyingxiang@buaa.edu.cn}{huyingxiang@buaa.edu.cn}}
	
	\vspace{3mm}
	\textsc{Institut f\"{u}r Diskrete Mathematik und Geometrie,\\ Technische Universit\"{a}t Wien, Wiedner Hauptstra\ss e 8-10,\\ 1040 Wien, Austria,}
	\email{\href{mailto:yingxiang.hu@tuwien.ac.at}{yingxiang.hu@tuwien.ac.at}}
	
	\vspace{3mm}
\textsc{Institut f\"{u}r Diskrete Mathematik und Geometrie,\\ Technische Universit\"{a}t Wien, Wiedner Hauptstra{\ss}e 8-10,\\ 1040 Wien, Austria,} \email{\href{mailto:mohammad.ivaki@tuwien.ac.at}{mohammad.ivaki@tuwien.ac.at}}

%	\textsc{Institut f\"{u}r Diskrete Mathematik und Geometrie,\\ Technische Universit\"{a}t Wien, Wiedner Hauptstra{\ss}e 8-10,\\ 1040 Wien, Austria,} \email{\href{mailto:mohammad.ivaki@tuwien.ac.at}{mohammad.ivaki@tuwien.ac.at}}
\end{document}